\documentclass[thmsa,12pt,sw20lart]{article}%
\usepackage{amsmath}
\usepackage{amsfonts}
\usepackage{amssymb}
\usepackage{graphicx}%
\newtheorem{theorem}{Theorem}[section]

\newtheorem{lemma}[theorem]{Lemma}

\newtheorem{remark}[theorem]{Remark}

\newenvironment{proof}[1][Proof]{\noindent\textbf{#1.} }{\ \rule{0.5em}{0.5em}}
\begin{document}

\author{Dong-Ho Tsai \ \ and\ \ \ Xiao-Liu Wang}
\title{On Length-preserving and Area-preserving Nonlocal Flow of Convex Closed\ Plane
Curves\thanks{AMS\ Subject Classifications: 53C44\ (primary),\ 35B40,\ 35K15,
35K55.}}
\maketitle

\begin{abstract}
For any $\alpha>0,$ we study $k^{\alpha}$-type length-preserving and
area-preserving nonlocal flow of convex closed plane curves and show that
these two types of flow evolve such curves into round circles in $C^{\infty}%
$-norm.$\ $Other relevant $k^{\alpha}$-type nonlocal flow is also discussed
when $\alpha\geq1.\ $

\end{abstract}

\section{Introduction.}

The curvature flow of plane curves, arising in many application fields, such
as phase transitions, image processing\ and smoothing,\ etc., has received a
lot of attention during the past several decades\ (see \cite{ST}\ for
example). In general, the evolution equation has the form:
\begin{equation}
\left\{
\begin{array}
[c]{l}%
\dfrac{\partial X}{\partial t}{\left(  u,t\right)  }=f\left(  k{\left(
u,t\right)  }\right)  \mathbf{N}_{in}{\left(  u,t\right)  ,}%
\vspace{3mm}%
\\
X{\left(  u,0\right)  }=X_{0}\left(  u\right)  ,\ \ \ u\in S^{1},
\end{array}
\right.  \label{eqn1.1}%
\end{equation}
where$\ X_{0}\left(  u\right)  \subset\mathbb{R}^{2}$ is a
given\ smooth\ closed curve, parametrized by $u\in S^{1},\ $and $X(u,t):S^{1}%
\times\lbrack0,T)\rightarrow{\mathbb{R}}^{2}\ $is a family of curves moving
along its inward normal direction$\ \mathbf{N}_{in}\left(  u,t\right)  $ with
given speed function $f\left(  k{(u,t)}\right)  ,\ $which is a strictly
increasing (parabolicity)\ function of the curvature$\ k{(u,t)}$ of
${X(u,t).\ }$

When $f\left(  k{(u,t)}\right)  =k(u,t)$, (\ref{eqn1.1}) is the famous
\emph{curve shortening flow}, which has been studied intensively by a great
number of authors for various conditions on the initial curve $X_{0}$.\ One
can see the book \cite{CZ} for literature. In particular, we mention the
papers \cite{AL, Ang, G, Ga1, GH}.\ When $f=\frac{1}{\alpha}k^{\alpha}%
\ $($\alpha\neq0\ $is a constant) and $X_{0}$ is a convex\footnote{Throughout
this paper, "convex"\ always means "uniformly convex",\ i.e., the curvature is
strictly positive everywhere.\ } simple closed curve, the flow (\ref{eqn1.1})
has been investigated in\ \cite{And1, And2}.\ Also in\ \cite{LPT1, LPT2, PT,
U1, U2}, they studied the case when$\ X_{0}$ is a locally convex\ non-simple
closed curve. For more about (\ref{eqn1.1}) and its counterpart in
high-dimensional space (more precisely, the \emph{mean curvature flow}) one
can find literature in the book \cite{Z}.\ 

Another class of interesting curvature flow is the so-called \emph{nonlocal
curvature flow}, the evolution equation of which takes the form:
\begin{equation}
\left\{
\begin{array}
[c]{l}%
{X}_{t}{\left(  u,t\right)  }=\left[  F(k{\left(  u,t\right)  })-\lambda
\left(  t\right)  \right]  \mathbf{N}_{in}{\left(  u,t\right)  },%
\vspace{3mm}%
\\
X{\left(  u,0\right)  }=X_{0}\left(  u\right)  ,\ \ \ u\in S^{1},
\end{array}
\right.  \label{eqn1.2}%
\end{equation}
where $X_{0}$ is a convex simple closed curve\ (we take its orientation to be
counterclockwise), $F(k)$ is a given function of the curvature satisfying the
parabolic condition$\ F^{\prime}(z)>0$ for all $z$ in its domain, and
$\lambda(t)$ is a function of time which depends on certain global
(nonlocal)\ quantities of $X\left(  \cdot,t\right)  $, say its\ length
$L(t),\ $enclosed area $A(t),$ or other possible global quantities like the
integral of curvature over the entire curve in certain ways.

The purpose of this paper is to study $k^{\alpha}$-type nonlocal flow
(\ref{eqn1.2})\ for convex closed curves with the speed function$\ F\left(
k\right)  -\lambda\left(  t\right)  $ given by $\ $
\begin{equation}
F\left(  k\right)  -\lambda\left(  t\right)  =k^{\alpha}-\dfrac{1}{{2\pi}}%
\int_{X\left(  \cdot,t\right)  }k^{\alpha+1}ds{,\ \ \ \alpha>0\ \ \ }%
\text{(LP)} \label{eqn1.3L}%
\end{equation}
or$\ $%
\begin{equation}
F\left(  k\right)  -\lambda\left(  t\right)  =k^{\alpha}-\dfrac{1}{{L}\left(
t\right)  }\int_{X\left(  \cdot,t\right)  }k^{\alpha}ds{,\ \ \ \alpha
>0\ \ \ }\text{(AP)}{.} \label{eqn1.3A}%
\end{equation}
Here $s$ is the arc\ length\ parameter of $X\left(  \cdot,t\right)  \ $and the
constant $\alpha>0$ is arbitrary. The abbreviation LP (AP)\ indicates that the
flow is\ "length-preserving" ("area-preserving"). We shall see shortly that
under the speed function\ (\ref{eqn1.3L}) the flow is length-preserving\ and
under (\ref{eqn1.3A}) the flow is area-preserving.\ Both flows have the
feature that it moves the point $p\in X\left(  \cdot,t\right)  \ $with maximal
curvature inward and the point $q\in X\left(  \cdot,t\right)  \ $with minimal
curvature outward.\ Moreover, they both\ decrease the\ isoperimetric ratio.

We note that the above\ two nonlocal flows are not suitable for\emph{
non-convex }simple closed\ curves due to the curvature term $k^{\alpha},$
since it may not be defined for $k<0.\ $Moreover, even if $k^{\alpha}%
,\ \alpha>0,$\ is defined for all $k\in(-\infty,\infty)\ $(say $\alpha
\in\mathbb{N}$), the nonlocal flow for a non-convex simple closed curve
$\gamma_{0}$ may easily develop\ self-intersections in short time. Thus we
have to confine to the convex case.\ 

When $\alpha=1,\ $(\ref{eqn1.3A}) is a gradient flow of the isoperimetric
deficit functional and it has been studied originally in \cite{Ga2}%
\ and\ recently in \cite{CLW}.\ They showed that the flow preserves
the\ convexity\ and enclosed area of $X_{0}$ and evolves it to a round circle
in $C^{\infty}$ sense as $t\rightarrow\infty.\ $As far as we know, Gage's
paper \cite{Ga2} seems to be the first\ one dealing with nonlocal flow of a
convex simple closed curve. Another very recent paper is \cite{MPW}, which
studied (\ref{eqn1.3A}) for the case when $\alpha\in\mathbb{N}$ is a positive integer.\ 

Also, when $\alpha=1$, (\ref{eqn1.3L}) has been studied in \cite{MZ}. They
showed that the flow preserves the\ convexity\ and\ length of$\ X_{0}$ and
evolves it to a round circle in $C^{\infty}$ sense as $t\rightarrow\infty.$
Another interesting nonlocal\ flow, with $F\left(  k\right)  =k,$ is due to
\cite{JP}, who studied a\ gradient flow of the isoperimetric ratio
functional.\ This\ gradient flow\ increases enclosed area\ and decreases
length of the evolving curve.\ One may say that, in some sense, it is the
fastest nonlocal flow (with $F\left(  k\right)  =k\ $and among all possible
choices of $\lambda\left(  t\right)  $)\ to evolve $X_{0}$ into a round circle.\ 

Even more recently there is another interesting paper \cite{BD}, which
discussed the AP nonlocal flow (with $F\left(  k\right)  =k$) for graphs (with
Dirichlet boundary condition).\ It has application to kinetic analysis of PET
(positron emission tomography)\ data.\ 

We first prove the following two $C^{\infty}\ $convergence results in Section
\ref{section2}.

\begin{theorem}
\label{thm-main-L}Assume$\ \alpha>0$\ and $X_{0}\left(  u\right)  ,\ u\in
S^{1},$ is a smooth convex closed curve\footnote{From now on, all closed
curves are assumed to be simple unless otherwise stated.}. Then under the LP
flow\ (\ref{eqn1.3L}), the flow exists\ and preserves length\ for all time
$t\in\lbrack0,\infty).$ Each$\ X\left(  \cdot,t\right)  $ remains smooth,
convex,\ and it converges to a round circle with radius $L\left(  0\right)
/2\pi\ $in $C^{\infty}\ $topology\ as $t\rightarrow\infty.$
\end{theorem}

\begin{theorem}
\label{thm-main-A}Assume$\ \alpha>0$\ and $X_{0}\left(  u\right)  ,\ u\in
S^{1},$ is a smooth convex closed curve. Then under the
AP\ flow\ (\ref{eqn1.3A}), the flow exists\ and preserves enclosed area\ for
all time $t\in\lbrack0,\infty).$ Each$\ X\left(  \cdot,t\right)  $ remains
smooth, convex,\ and it converges to a round circle with radius $\sqrt
{A\left(  0\right)  /\pi}\ $in $C^{\infty}\ $topology\ as $t\rightarrow
\infty.$
\end{theorem}

Later in Section \ref{general} we also consider other relevant $k^{\alpha}%
$-type nonlocal flow for$\ \alpha\geq1\ $(see Theorem \ref{thm-2}). This
flow,\ which is similar to \cite{JP}, has the property that it is
area-increasing and\ length-decreasing.\ The methods of proof for Theorem
\ref{thm-main-L}\ and Theorem \ref{thm-main-A}\ can be applied to Theorem
\ref{thm-2} without essential changes.

One can also study $1/k^{\alpha}$-type (i.e., $F\left(  k\right)
=-1/k^{\alpha}$) nonlocal flow for $\alpha>0.\ $We will not discuss it
here.\ For the case when $F(k)=-k^{-1},$ we mention the papers \cite{MC, PY,
PZ}.\ Also see \cite{LT} for other types of nonlocal flow.\ 

Finally, we point out that a lot of the methods of proof used in the
well-known curve shortening flow (or other related local curvature flows)\ are
not quite applicable to the nonlocal flows we plan to discuss here. Thus we
need to use different\ techniques. \ 

There are also many interesting nonlocal flows in the higher dimensional case,
we mention the papers \cite{And3, ES, H, M1, M2}.

\section{$k^{\alpha}$-type LP and AP flow; the proof of Theorem
\ref{thm-main-L}\ and Theorem \ref{thm-main-A}.\label{section2}}

\subsection{Short time existence of the flow (\ref{eqn1.3L}) (or
(\ref{eqn1.3A})).\label{existence}}

The nonlocal flow (\ref{eqn1.3L}) (or (\ref{eqn1.3A})) is nonlinear in
nature.\ It is fully-nonlinear if we look at the evolution of the support
function and quasilinear if we look at the evolution of the curvature.\ Hence
one may have concern about the short time existence of a solution to these two
flows.\ However, in view of the recent various nonlocal flow papers in the
literature, it is not a problem.

More precisely, since we assume that the initial curve $X_{0}\left(  u\right)
\ $is smooth and convex (with curvature strictly positive everywhere), the
equation\ is \emph{strictly parabolic initially}, and the nonlinear nature of
the equation is not really an obstacle --\ for example, one can look at the
treatment of short time existence in Chapter 3 of \cite{B}.\ The nonlocal term
does add a little complication, but it is of a very mild form and can be
easily handled.\ See the argument in Section 7 of \cite{M2}. Another method is
to use the classical Leray-Schauder fixed point theory.\ See \cite{MPW} for details.\ 

Finally, we note that the solution to the flow (\ref{eqn1.3L}) (or
(\ref{eqn1.3A})) is also unique.\ Again, this is due to the fact that the
equation is strictly parabolic initially or one can follow similar treatment
as in \cite{MPW}.\ Also see \cite{ES, M1, M2}.

Thus we conclude that there is a unique smooth\ solution to the flow
(\ref{eqn1.3L}) (or (\ref{eqn1.3A})) on $S^{1}\times\lbrack0,T)$ for some
short time $T>0.\ $Moreover,\ as $X_{0}\left(  u\right)  \ $is convex, by
continuity,\ we may assume that each $X\left(  \cdot,t\right)  \ $is also
convex for $t\in\lbrack0,T).\ $Later, we shall show that the convexity is
preserved as long as the solution to the flow exists. $\ $

\subsection{Upper bound of the curvature;\ Tso's method.}

In this section, we will use Tso's method \cite{T}\ to obtain a time-dependent
upper bound of the curvature for $\alpha>0$.

By Section \ref{existence}, we have a unique smooth convex solution $X\left(
u,t\right)  \ $to both nonlocal flows on $S^{1}\times\lbrack0,T)$.\ According
to \cite{And1, Ang, GH} and many other papers, we can use the outward normal
angle $\theta\in S^{1}=\left[  0,2\pi\right]  $ to parameterize $X\left(
\cdot,t\right)  $.\ By $X\left(  \theta,t\right)  $ we mean the unique point
on $X\left(  \cdot,t\right)  $ at which its outward normal $\mathbf{N}%
_{out}\mathbf{\ }$is\ given\ by$\ \left(  \cos\theta,\sin\theta\right)
.$\ With this, all evolution equations of the flow\ can be expressed in
$\left(  \theta,t\right)  \ $coordinates\ and the evolution of the curvature
$k\left(  \theta,t\right)  \ $of $X\left(  \theta,t\right)  $ is given
by\ (see\ \cite{GH, Ga2} for computation)
\begin{equation}
\left\{
\begin{array}
[c]{l}%
k_{t}\left(  \theta,t\right)  =k^{2}\left(  \theta,t\right)  \left[  \left(
k^{\alpha}\right)  _{\theta\theta}\left(  \theta,t\right)  +k^{\alpha}\left(
\theta,t\right)  -\lambda\left(  t\right)  \right]  ,%
\vspace{3mm}%
\\
k(\theta,0)=k_{0}\left(  \theta\right)  >0,\ \ \ \left(  \theta,t\right)  \in
S^{1}\times\lbrack0,T),
\end{array}
\right.  \label{dkdt}%
\end{equation}
where $k_{0}\left(  \theta\right)  >0$ is the initial curvature of
$X_{0}\left(  \theta\right)  \ $and%
\begin{equation}
\lambda\left(  t\right)  =\dfrac{1}{{2\pi}}\int_{X(\cdot,\,t)}k^{\alpha
+1}ds=\dfrac{1}{{2\pi}}%
{\displaystyle\int_{0}^{2\pi}}
k^{\alpha}\left(  \theta,t\right)  d\theta{\ \ \ }\text{(LP)} \label{lamda-L}%
\end{equation}
or%
\begin{equation}
\lambda\left(  t\right)  =\dfrac{1}{{L}\left(  t\right)  }\int_{X\left(
\cdot,t\right)  }k^{\alpha}ds=\left(
{\displaystyle\int_{0}^{2\pi}}
\frac{1}{k\left(  \theta,t\right)  }d\theta\right)  ^{-1}%
{\displaystyle\int_{0}^{2\pi}}
k^{\alpha-1}\left(  \theta,t\right)  d\theta{\ \ \ }\text{(AP).}
\label{lamda-A}%
\end{equation}
Since $X_{0}\left(  \theta\right)  $ is a closed curve, $k_{0}\left(
\theta\right)  $ satisfies the integral condition%
\begin{equation}
\int_{0}^{2\pi}\frac{\cos\theta}{k_{0}(\theta)}d\theta=\int_{0}^{2\pi}%
\frac{\sin\theta}{k_{0}(\theta)}d\theta=0. \label{int}%
\end{equation}
By (\ref{dkdt}), one can easily check that (\ref{int}) is preserved under both
flows. \ 

One interesting comparison is in order. In the LP flow $\lambda\left(
t\right)  $ is the average of $k^{\alpha}$ over $S^{1}\ $%
with\ respect\ to\ the angle$\ \theta\in\left[  0,2\pi\right]  $, while in the
AP flow\ $\lambda\left(  t\right)  $ is the average of $k^{\alpha}$ over
$X\left(  \cdot,t\right)  \ $with respect to the arc\ length $s\in\left[
0,L\left(  t\right)  \right]  $.\ Also note that, for a convex closed curve
$\gamma\subset\mathbb{R}^{2}\ $with\ curvature\ $k,\ $we have%
\begin{equation}
\dfrac{1}{{L}}\int_{\gamma}k^{\alpha}ds\leq\dfrac{1}{{2\pi}}\int_{\gamma
}k^{\alpha+1}ds,\ \ \ {\alpha>0} \label{kkL}%
\end{equation}
due to the H\"{o}lder inequality.\ Thus, for the same curve, the nonlocal term
in\ the LP\ flow has\textbf{ }more outward effect than that in\ the AP\ flow.

Similar to Theorem 4.1.4 of \cite{GH}, one can verify that the
nonlocal\textbf{ }flow (\ref{eqn1.3L}) is equivalent\ to the equation
(\ref{dkdt}) with $\lambda\left(  t\right)  $ given by\ (\ref{lamda-L});
and\ the nonlocal flow (\ref{eqn1.3A}) is equivalent\ to the equation
(\ref{dkdt}) with\ $\lambda\left(  t\right)  \ $given by\ (\ref{lamda-A}%
).\ \ Both have\ smooth$\ $initial$\ $data$\ k_{0}\left(  \theta\right)  >0$
satisfying (\ref{int}). In view of this, from now on,\ we can just focus on
the curvature equation (\ref{dkdt}) and study its long time behavior.\ This
has become a rather standard way of studying curvature flows.\ 

By Lemma 2.2 in \cite{Ga2},\ for the LP case,\ we have
\begin{equation}
\frac{dL}{dt}\left(  t\right)  =0,\ \ \ \frac{dA}{dt}\left(  t\right)
=-\int_{0}^{2\pi}k^{\alpha-1}\left(  \theta,t\right)  d\theta+\dfrac
{{L\!}\left(  t\right)  }{2\pi}\int_{0}^{2\pi}k^{\alpha}\left(  \theta
,t\right)  d\theta\geq0,\ \ \ t\in\lbrack0,T), \label{LA1}%
\end{equation}
and for the AP case,\ we have%
\begin{equation}
\frac{dA}{dt}\left(  t\right)  =0,\ \ \ \frac{dL}{dt}\left(  t\right)
=-\int_{0}^{2\pi}k^{\alpha}\left(  \theta,t\right)  d\theta+\dfrac{2\pi}%
{{L}\left(  t\right)  }\int_{0}^{2\pi}k^{\alpha-1}\left(  \theta,t\right)
d\theta\leq0,\ \ \ t\in\lbrack0,T). \label{LA2}%
\end{equation}
For $\alpha>0,\ $the inequalities in (\ref{LA1})\ and (\ref{LA2})\ are both
due to (\ref{kkL}) and the identity $ds=k^{-1}d\theta.$

By (\ref{LA1})\ and (\ref{LA2}),\ in both flows, the evolving curve $X\left(
\cdot,t\right)  \ $is becoming more and more circular in the sense that the
isoperimetric ratio $L^{2}\left(  t\right)  /4\pi A\left(  t\right)  $
is$\ $decreasing in time $t\in\lbrack0,T).$ To show $C^{\infty}$ convergence
of $X\left(  \cdot,t\right)  $, we need to do curvature estimate, which is the
main work of the paper.\ 

To go further, we need the following \emph{Bonnesen inequality} for convex
closed plane\ curves (see the book by \cite{S},\ p.\ 324): For any\ convex
closed curve\textbf{ }$\gamma$ in the plane,\ there holds the inequality%
\begin{equation}
\rho L-A-\pi\rho^{2}\geq0\ \ \ \text{for all\ \ }r\leq\rho\leq R, \label{BONN}%
\end{equation}
where $r\ $and$\ R\ $are the\ \emph{inradius} and \emph{outradius}
(\emph{circumradius})\ of $\gamma\ $respectively.

By (\ref{BONN}), we have
\begin{equation}
0<\frac{L-\sqrt{L^{2}-4\pi A}}{2\pi}\leq r\leq R\leq\frac{L+\sqrt{L^{2}-4\pi
A}}{2\pi},
\end{equation}
and therefore obtain the inequality
\begin{equation}
1\leq\frac{R}{r}\leq\frac{L+\sqrt{L^{2}-4\pi A}}{L-\sqrt{L^{2}-4\pi A}%
}=\left(  \sqrt{I}+\sqrt{I-1}\right)  ^{2},\ \ \ \text{where\ \ \ }%
I=\frac{L^{2}}{4\pi A}\geq1. \label{bonn}%
\end{equation}
\ \ \ \ 

We now can use Tso's method to prove a time-dependent upper bound of the curvature.

\begin{lemma}
\label{lem-1}Under\ the\ LP flow\ (\ref{eqn1.3L}) on $S^{1}\times
\lbrack0,T)\ $with\ $\alpha>0,\ $there exists a constant $C\left(  T\right)
>0\ $depending on the time $T\ $such that%
\begin{equation}
0<v\left(  \theta,t\right)  \leq C\left(  T\right)  ,\text{ \ \ }\forall\text{
}\left(  \theta,t\right)  \in S^{1}\times\lbrack0,T)
\end{equation}
i.e.,%
\begin{equation}
0<k\left(  \theta,t\right)  \leq C\left(  T\right)  ,\ \ \ \forall\ \left(
\theta,t\right)  \in S^{1}\times\lbrack0,T). \label{k-upper-bd}%
\end{equation}
Here $v\left(  \theta,t\right)  =k^{\alpha}\left(  \theta,t\right)  $ and
$k\left(  \theta,t\right)  $ is the curvature of $X\left(  \cdot,t\right)  .$
The same result holds for the AP\ flow\ (\ref{eqn1.3A}) on $S^{1}\times
\lbrack0,T)\ $with\ $\alpha>0.$
\end{lemma}

%

\proof
The evolution of $v\left(  \theta,t\right)  =k^{\alpha}\left(  \theta
,t\right)  $ is given by%
\begin{equation}
v_{t}=\alpha v^{p}\left(  v_{\theta\theta}+v-\lambda\left(  t\right)  \right)
,\ \ \ p=1+\frac{1}{\alpha}>1,\ \ \ \left(  \theta,t\right)  \in S^{1}%
\times\lbrack0,T), \label{v-eq}%
\end{equation}
where$\ \lambda\left(  t\right)  >0\ $is from (\ref{lamda-L})
or\ (\ref{lamda-A}). Since the isoperimetric ratio $I\left(  t\right)
=L^{2}\left(  t\right)  /4\pi A\left(  t\right)  \ $of $X\left(
\cdot,t\right)  \ $is decreasing in time in both flows, we have$\ 1\leq
I\left(  t\right)  \leq I\left(  0\right)  $\ for all $t\in\lbrack0,T).\ $By
the Bonnesen inequality (\ref{bonn}), we have%
\begin{equation}
1\leq\frac{r_{o}\left(  t\right)  }{r_{i}\left(  t\right)  }\leq\left(
\sqrt{I\left(  t\right)  }+\sqrt{I\left(  t\right)  -1}\right)  ^{2}%
\leq\left(  \sqrt{I\left(  0\right)  }+\sqrt{I\left(  0\right)  -1}\right)
^{2}:=\sigma,\ \ \ t\in\lbrack0,T) \label{Bonn}%
\end{equation}
where $r_{o}\left(  t\right)  \ $and $r_{i}\left(  t\right)  $ are the
outradius and inradius of $X\left(  \cdot,t\right)  $ respectively.\ Hence
$r_{o}\left(  t\right)  \leq\sigma r_{i}\left(  t\right)  \ $and by $A\left(
0\right)  \leq A\left(  t\right)  \leq\pi r_{o}^{2}\left(  t\right)  ,\ $we
have\ $r_{i}\left(  t\right)  \geq\sigma^{-1}\sqrt{A\left(  0\right)  /\pi}%
$\ for all $t\in\lbrack0,T),$\ i.e.,$\ $the inradius has
time-independent\ positive lower bound.

Let $E\left(  0\right)  $ be a circle enclosed by $X_{0}.\ $By the maximum
principle, if we shrink $E\left(  0\right)  $ by the $k^{\alpha}$-contraction
flow\ (without the nonlocal term $\lambda\left(  t\right)  $), then $E\left(
t\right)  $ is enclosed by ${X(\cdot,t)}$ for all $t\in\lbrack0,T).$ The
radius $r\left(  t\right)  \ $of $E\left(  t\right)  $ is given by%
\begin{equation}
r\left(  t\right)  =\left[  r^{1+\alpha}\left(  0\right)  -\left(
1+\alpha\right)  t\right]  ^{1/\left(  1+\alpha\right)  },\ \ \ t\in
\lbrack0,r^{1+\alpha}\left(  0\right)  /\left(  1+\alpha\right)  ). \label{rt}%
\end{equation}
Now we take$\ E\left(  0\right)  $ to be an inscribed circle of $X_{0}\left(
\cdot\right)  \ $with inradius $r\left(  0\right)  =r_{i}\left(  0\right)
.\ $If we use the center of $E\left(  0\right)  $ as the origin\ $O$, then the
support function$\ u\left(  \theta,t\right)  $ (with respect to $O$) of
$X\left(  \cdot,t\right)  \ $is defined on $S^{1}\times\lbrack0,T)\ $and by
the above we have\ (since $X\left(  \cdot,t\right)  $ encloses $E\left(
t\right)  $)%
\begin{equation}
u\left(  \theta,t\right)  \geq\left[  r^{1+\alpha}\left(  0\right)  -\left(
1+\alpha\right)  t\right]  ^{1/\left(  1+\alpha\right)  },\ \ \ \left(
\theta,t\right)  \in S^{1}\times\lbrack0,\min\left\{  r^{1+\alpha}\left(
0\right)  /\left(  1+\alpha\right)  ,T\right\}  ), \label{ur}%
\end{equation}
where now $r\left(  0\right)  =r_{i}\left(  0\right)  \geq\sigma^{-1}%
\sqrt{A\left(  0\right)  /\pi}$.

Since, in general, $T>0$ is a small time,\ we may assume that$\ $%
\begin{equation}
T\leq\frac{1}{2+2\alpha}\left(  \frac{1}{\sigma}\sqrt{\frac{A\left(  0\right)
}{\pi}}\right)  ^{1+\alpha}:=T_{1} \label{T}%
\end{equation}
(the case $T>T_{1}\ $can be handled similarly).\ \ By (\ref{ur}%
)\ and$\ r\left(  0\right)  \geq\sigma^{-1}\sqrt{A\left(  0\right)  /\pi}$, we
have$\ u\left(  \theta,t\right)  \geq2\beta\ $on$\ S^{1}\times\lbrack
0,T),\ $where $\beta>0\ $is a constant\ depending only on the initial curve
$X_{0}\left(  \cdot\right)  ,$ given by
\begin{equation}
\beta=\left(  \frac{1}{2}\right)  ^{\left(  2+\alpha\right)  /\left(
1+\alpha\right)  }\frac{1}{\sigma}\sqrt{\frac{A\left(  0\right)  }{\pi}%
},\ \ \ \sigma=\left(  \sqrt{I\left(  0\right)  }+\sqrt{I\left(  0\right)
-1}\right)  ^{2}. \label{beta}%
\end{equation}
Moreover, as the length $L\left(  t\right)  \ $of $X\left(  \cdot,t\right)  $
is decreasing in both flows, there is a constant $C>0$ depending only on
$X_{0}\left(  \cdot\right)  $ such that$\ u\left(  \theta,t\right)  \leq
C\ $on$\ S^{1}\times\lbrack0,T).\ $Hence
\begin{equation}
0<2\beta\leq u\left(  \theta,t\right)  \leq C\ \ \ \text{on}\ \ \ S^{1}%
\times\lbrack0,T). \label{uC}%
\end{equation}
Note that both $\beta\ $and$\ C\ $are\ independent\ of time.\ 

Following \cite{T}, we consider the evolution of the quantity$\ $
\begin{equation}
Q\left(  \theta,t\right)  =\frac{v\left(  \theta,t\right)  }{u\left(
\theta,t\right)  -\beta},\ \ \ v\left(  \theta,t\right)  =k^{\alpha}\left(
\theta,t\right)  ,\ \ \ \left(  \theta,t\right)  \in S^{1}\times\lbrack0,T).
\label{Q}%
\end{equation}
Under the nonlocal\ flow, the evolution of the support function is%
\begin{equation}
u_{t}\left(  \theta,t\right)  =-v\left(  \theta,t\right)  +\lambda\left(
t\right)  =-k^{\alpha}\left(  \theta,t\right)  +\lambda\left(  t\right)
,\ \ \ \left(  \theta,t\right)  \in S^{1}\times\lbrack0,T), \label{supp}%
\end{equation}
where$\ \lambda\left(  t\right)  >0\ $is from (\ref{lamda-L})
or\ (\ref{lamda-A}).\ We compute the evolution of $Q\ $on\ $S^{1}\times
\lbrack0,T)\ $to get%
\[
Q_{t}=\frac{v\left(  v-\lambda\left(  t\right)  \right)  }{\left(
u-\beta\right)  ^{2}}+\frac{\alpha v^{p}\left(  v_{\theta\theta}%
+v-\lambda\left(  t\right)  \right)  }{u-\beta},\ \ \ p=1+\frac{1}{\alpha}%
\]
and by%
\[
Q_{\theta}=\frac{v_{\theta}}{u-\beta}-\frac{u_{\theta}v}{\left(
u-\beta\right)  ^{2}},\ \ \ Q_{\theta\theta}=\frac{v_{\theta\theta}}{u-\beta
}-\frac{2u_{\theta}v_{\theta}}{\left(  u-\beta\right)  ^{2}}+\left(
\frac{2u_{\theta}^{2}}{\left(  u-\beta\right)  ^{3}}-\frac{u_{\theta\theta}%
}{\left(  u-\beta\right)  ^{2}}\right)  v
\]
we conclude%
\begin{align}
Q_{t}  &  =Q^{2}-\frac{\lambda\left(  t\right)  v}{\left(  u-\beta\right)
^{2}}+\alpha v^{p}\left[  Q_{\theta\theta}+\frac{2u_{\theta}v_{\theta}%
}{\left(  u-\beta\right)  ^{2}}-\left(  \frac{2u_{\theta}^{2}}{\left(
u-\beta\right)  ^{3}}-\frac{u_{\theta\theta}}{\left(  u-\beta\right)  ^{2}%
}\right)  v+Q-\frac{\lambda\left(  t\right)  }{u-\beta}\right] \nonumber\\
&  \leq\left\{
\begin{array}
[c]{l}%
\alpha v^{p}Q_{\theta\theta}+\alpha v^{p}\frac{2u_{\theta}}{u-\beta}\left(
Q_{\theta}+\frac{u_{\theta}v}{\left(  u-\beta\right)  ^{2}}\right)
\vspace{3mm}%
\\
-\alpha v^{p+1}\left(  \frac{2u_{\theta}^{2}}{\left(  u-\beta\right)  ^{3}%
}-\frac{1}{k}\frac{1}{\left(  u-\beta\right)  ^{2}}+\frac{u}{\left(
u-\beta\right)  ^{2}}\right)  +\left(  \alpha v^{p}+Q\right)  Q
\end{array}
\right. \nonumber\\
&  =\alpha v^{p}Q_{\theta\theta}+\alpha v^{p}\frac{2u_{\theta}}{u-\beta
}Q_{\theta}+Q^{2}\left[  \left(  \alpha+1\right)  -\alpha\beta\left(
u-\beta\right)  ^{1/\alpha}Q^{1/\alpha}\right] \nonumber\\
&  \leq\alpha v^{p}Q_{\theta\theta}+\alpha v^{p}\frac{2u_{\theta}}{u-\beta
}Q_{\theta}+Q^{2}\left[  \left(  \alpha+1\right)  -\alpha\beta^{1+1/\alpha
}Q^{1/\alpha}\right]  , \label{Qt}%
\end{align}
where we have used the familiar identity $u_{\theta\theta}+u=1/k\ $and the
inequality$\ u-\beta\geq\beta\ $on\ $S^{1}\times\lbrack0,T)$ in the
above.\ Also note that in the above estimate\ we have thrown away the two
favorable terms containing $\lambda\left(  t\right)  \ $due\ to$\ $the
negative sign in front of them.

Let$\ Q_{0}\left(  \alpha,\beta\right)  =\left[  2\left(  \alpha+1\right)
/\left(  \alpha\beta^{1+1/\alpha}\right)  \right]  ^{\alpha}>0.\ $It is a
constant depending only on $\alpha,\ \beta.$\ Then\ whenever%
\[
Q_{\max}\left(  t\right)  =Q\left(  \theta\left(  t\right)  ,t\right)  \geq
Q_{0}\left(  \alpha,\beta\right)  ,\ \ \ t\in\lbrack0,T),
\]
we have
\[
\frac{d}{dt}Q_{\max}\left(  t\right)  \leq Q_{t}\left(  \theta\left(
t\right)  ,t\right)  \leq-\left(  \alpha+1\right)  Q_{\max}^{2}\left(
t\right)  ,\ \ \ t\in\lbrack0,T).
\]
Comparison principle implies
\begin{equation}
Q_{\max}\left(  t\right)  \leq\max\left\{  Q_{0}\left(  \alpha,\beta\right)
,\ \ \frac{1}{\left(  \alpha+1\right)  t}\right\}  ,\ \ \ \forall\ t\in\left(
0,T\right)  . \label{Qmax}%
\end{equation}
In particular, we get%
\[
\frac{k^{\alpha}\left(  \theta,t\right)  }{u\left(  \theta,t\right)  -\beta
}\leq\max\left\{  Q_{0}\left(  \alpha,\beta\right)  ,\ \;\frac{1}{\left(
\alpha+1\right)  t}\right\}  ,\ \ \ \forall\ \left(  \theta,t\right)  \in
S^{1}\times\left(  0,T\right)
\]
and then by (\ref{uC}) we conclude
\begin{equation}
k_{\max}^{\alpha}\left(  t\right)  \leq\max\left\{  \left(  C-\beta\right)
Q_{0}\left(  \alpha,\beta\right)  ,\ \;\frac{C-\beta}{\left(  \alpha+1\right)
t}\right\}  ,\ \ \ \forall\ t\in\left(  0,T\right)  , \label{kmax}%
\end{equation}
which implies that near time $T\ $the curvature is bounded above by a constant
$C\left(  T\right)  \ $depending on $T,\ $i.e., the curvature will not
blow\ up as $t\rightarrow T.\ $The proof is done.$%
\hfill
\square$

\subsection{Lower bound of the curvature.}

We first observe the following:\ 

\begin{lemma}
\label{lem-2}Under the LP flow (\ref{eqn1.3L}) on $S^{1}\times\lbrack
0,T)\ $with $\alpha>0,$ there holds the estimate%
\begin{equation}
\max_{S^{1}\times\left[  0,t\right]  }\Psi\leq\max\left\{  \max_{S^{1}%
\times\left[  0,t\right]  }v^{2},\ \ \max_{S^{1}\times\left\{  0\right\}
}\Psi\right\}  ,\ \ \ \forall\ t\in\lbrack0,T), \label{v-grad-Ben}%
\end{equation}
where$\ v=k^{\alpha}\ $and$\ \Psi=v^{2}+v_{\theta}^{2}.$ In particular, we
have%
\begin{equation}
\left\vert v_{\theta}\left(  \theta,t\right)  \right\vert \leq C\left(
T\right)  ,\ \ \ \forall\ \left(  \theta,t\right)  \in S^{1}\times\lbrack0,T).
\label{v-grad-est}%
\end{equation}
The same result holds for the AP\ flow\ (\ref{eqn1.3A}) on $S^{1}\times
\lbrack0,T)$ with $\alpha>0.$
\end{lemma}

%

\proof
The proof of (\ref{v-grad-Ben})\ is analogous to that in Lemma I1.12 in
\cite{And1}. The evolution of $\Psi=v^{2}+v_{\theta}^{2}\ $is%
\begin{align}
\partial_{t}\Psi &  =2\alpha v^{p+1}\left(  v_{\theta\theta}+v-\lambda\left(
t\right)  \right)  +\left[  \alpha v^{p}\Psi_{\theta\theta}-2\alpha
v^{p}v_{\theta\theta}\left(  v_{\theta\theta}+v\right)  \right]  +\alpha
pv^{p-1}v_{\theta}\Psi_{\theta}-2\alpha pv^{p-1}v_{\theta}^{2}\lambda\left(
t\right)  \nonumber\\
&  \leq\alpha v^{p}\Psi_{\theta\theta}+\alpha pv^{p-1}v_{\theta}\Psi_{\theta
}+2\alpha v^{p+1}\left(  v_{\theta\theta}+v\right)  -2\alpha v^{p}%
v_{\theta\theta}\left(  v_{\theta\theta}+v\right)  ,\text{ \ \ }p=1+\frac
{1}{\alpha}>1\label{evo}%
\end{align}
where in the above inequality we have thrown away the terms containing
$\lambda\left(  t\right)  .$ On the time interval $\left[  0,t\right]  ,$ let
$\sigma\ $be the constant$\ \max_{S^{1}\times\left[  0,t\right]  }v^{2}%
.\ $Then whenever $\Psi_{\max}\left(  s\right)  =\Psi\left(  \theta
_{s},s\right)  >\sigma$ at any time $s\in\left[  0,t\right]  ,$ we have\
\[
\Psi_{\theta}\left(  \theta_{s},s\right)  =2v_{\theta}\left(  \theta
_{s},s\right)  \left(  v_{\theta\theta}\left(  \theta_{s},s\right)  +v\left(
\theta_{s},s\right)  \right)  =0,\ \ \ \ \ \Psi_{\theta\theta}\left(
\theta_{s},s\right)  \leq0
\]
and since $\Psi\left(  \theta_{s},s\right)  >\sigma,$ we must have\ $v_{\theta
}\left(  \theta_{s},s\right)  \neq0$ and so $\left(  v_{\theta\theta}\left(
\theta_{s},s\right)  +v\left(  \theta_{s},s\right)  \right)  =0.\ $Now by the
evolution inequality (\ref{evo}) we get$\ \partial_{t}\Psi\left(  \theta
_{s},s\right)  \leq0\ $at$\ \left(  \theta_{s},s\right)  .\ $As a consequence
of the maximum principle, we have estimate (\ref{v-grad-Ben}). Finally, the
estimate$\ $(\ref{v-grad-est}) follows from (\ref{v-grad-Ben})\ and
(\ref{k-upper-bd}).$%
\hfill
\square$

\ 

To obtain a time-dependent lower bound of the curvature for $\alpha>0$, we use
the idea from p. 64 of \cite{MZ} (which is for the case $\alpha=1$).\ \ We have:

\begin{lemma}
\label{lem-3}Under the LP flow (\ref{eqn1.3L}) on $S^{1}\times\lbrack
0,T)$\ with $\alpha>0,$ there exists a constant $c\left(  T\right)  >0$ such
that
\begin{equation}
k\left(  \theta,t\right)  \geq c\left(  T\right)  >0,\ \ \ \forall
\text{\ }\left(  \theta,t\right)  \in S^{1}\times\lbrack0,T).
\label{k-lower-bd}%
\end{equation}
The same result holds for the AP\ flow\ (\ref{eqn1.3A}) on $S^{1}\times
\lbrack0,T)$\ with $\alpha>0.$
\end{lemma}

%

\proof
Consider the quantity
\begin{equation}
\Phi\left(  \theta,t\right)  =\frac{1}{k\left(  \theta,t\right)  }%
-\frac{L\left(  t\right)  }{2\pi}-\frac{1}{2\pi}\int_{0}^{t}\int_{0}^{2\pi
}k^{\alpha}\left(  \theta,\sigma\right)  d\theta d\sigma,\ \ \ \left(
\theta,t\right)  \in S^{1}\times\lbrack0,T),
\end{equation}
where $\max_{\theta\in S^{1}}\Phi\left(  \theta,0\right)  \geq0.\ $In both
flows,\ we have
\begin{align}
\Phi_{t}  &  =-\left(  k^{\alpha}\right)  _{\theta\theta}-k^{\alpha}%
+\lambda\left(  t\right)  -\left(  \lambda\left(  t\right)  -\frac{1}{2\pi
}\int_{0}^{2\pi}k^{\alpha}d\theta\right)  -\frac{1}{2\pi}\int_{0}^{2\pi
}k^{\alpha}\left(  \theta,t\right)  d\theta\nonumber\\
&  =\alpha k^{\alpha+1}\left(  \theta,t\right)  \Phi_{\theta\theta}\left(
\theta,t\right)  -\alpha(\alpha+1)k^{\alpha+2}\left(  \theta,t\right)
\Phi_{\theta}^{2}\left(  \theta,t\right)  -k^{\alpha}\left(  \theta,t\right)
\leq\alpha k^{\alpha+1}\Phi_{\theta\theta} \label{phi}%
\end{align}
and the maximum principle implies$\ $ \ $\ $%
\begin{equation}
\frac{1}{k\left(  \theta,t\right)  }\leq\max_{\theta\in S^{1}}\left(  \frac
{1}{k_{0}\left(  \theta\right)  }\right)  +\frac{L\left(  t\right)  -L\left(
0\right)  }{2\pi}+\frac{1}{2\pi}\int_{0}^{t}\int_{0}^{2\pi}k^{\alpha}\left(
\theta,\sigma\right)  d\theta d\sigma\label{1/k-2}%
\end{equation}
for all $\left(  \theta,t\right)  \in S^{1}\times\lbrack0,T).\ $By
(\ref{k-upper-bd}) and the fact that $L\left(  t\right)  $ is decreasing, we
get estimate (\ref{k-lower-bd}).$%
\hfill
\square$

\subsection{Long time existence and $C^{\infty}\ $convergence of the flow.}

The curvature estimates established so far and the parabolic regularity theory
implies the following:

\begin{theorem}
\label{thm-1}(\emph{Long time existence of the flow.})\ For any $\alpha
>0,\ $the LP flow (\ref{eqn1.3L}) has a smooth\ convex solution\ $X\left(
\cdot,t\right)  $ defined for all time $t\in\lbrack0,\infty).$
Moreover,\ there exists a constant $C>0,\ $independent\ of time,\ such that%
\begin{equation}
0<v\left(  \theta,t\right)  \leq C\ \ \ \text{and\ \ \ }\left\vert v_{\theta
}(\theta,t)\right\vert \leq C,\text{ \ \ }\forall\ \left(  \theta,t\right)
\in S^{1}\times\lbrack0,\infty), \label{kv}%
\end{equation}
i.e.,%
\begin{equation}
0<k\left(  \theta,t\right)  \leq C,\ \ \ \forall\ \left(  \theta,t\right)  \in
S^{1}\times\lbrack0,\infty), \label{k-up-time-ind}%
\end{equation}
where\ $v=k^{\alpha}.\ $The same result holds for the AP\ flow\ (\ref{eqn1.3A}%
)\ with $\alpha>0.$
\end{theorem}

%

\proof
On the domain$\ S^{1}\times\lbrack0,T),\ $the curvature $k\left(
\theta,t\right)  \ $has uniform positive upper and lower bounds.\ Parabolic
regularity then implies that all space-time derivatives of $k\left(
\theta,t\right)  $ remain bounded on the domain$\ S^{1}\times\lbrack
0,T).\ $Therefore, the evolving curve has a smooth convex\ limit\ $X\left(
\cdot,T\right)  $ as $t\rightarrow T$ and one can use$\ X\left(
\cdot,T\right)  $ as a new initial curve to continue the flow.\ Hence the flow
is defined on $S^{1}\times\lbrack0,\infty).\ $It remains to explain
(\ref{kv}).$\ $

We now have $T=\infty$ in Lemma \ref{lem-1}.\ For any$\ \xi\geq0$ we can
apply\ Tso's argument on the interval $\left[  \xi,\xi+T_{1}\right]
\ $($T_{1}$ is from (\ref{T}), which is independent\ of $\xi$).\ The support
function $u\ $is now with respect to\ the center of an inscribed circle of
$X\left(  \cdot,\xi\right)  \ $with radius $r_{i}\left(  \xi\right)
\geq\sigma^{-1}\sqrt{A\left(  0\right)  /\pi}$. The constants$\ C\ $%
and$\ \beta\ $are the same as before (both are independent\ of $\xi$)$\ $and
we still have inequality\ (\ref{uC})\ on\ $S^{1}\times\left[  \xi,\xi
+T_{1}\right]  .\ $Since the time interval is $\left[  \xi,\xi+T_{1}\right]
$, the estimate (\ref{Qmax}) now becomes%
\begin{equation}
Q_{\max}\left(  t\right)  \leq\max\left\{  Q_{0}\left(  \alpha,\beta\right)
,\ \ \frac{1}{\left(  \alpha+1\right)  \left(  t-\xi\right)  }\right\}
,\ \ \ \forall\ t\in\left(  \xi,\xi+T_{1}\right)
\end{equation}
and we get $\ \ $%
\[
0<k_{\max}^{\alpha}\left(  t\right)  \leq\max\left\{  \left(  C-\beta\right)
Q_{0}\left(  \alpha,\beta\right)  ,\ \;\frac{C-\beta}{\left(  \alpha+1\right)
\left(  t-\xi\right)  }\right\}  ,\ \ \ \forall\ t\in\left(  \xi,\xi
+T_{1}\right)  .
\]
In particular, we get%
\[
0<k_{\max}^{\alpha}\left(  \xi+T_{1}\right)  \leq\max\left\{  \left(
C-\beta\right)  Q_{0}\left(  \alpha,\beta\right)  ,\ \;\frac{C-\beta}{\left(
\alpha+1\right)  T_{1}}\right\}  .
\]
As $\xi\geq0$ is arbitrary, the first estimate in (\ref{kv}) is
confirmed.\ The second estimate follows from (\ref{v-grad-Ben})\ in Lemma
\ref{lem-2}.$%
\hfill
\square\ $

\ \ \ \ 

Until now, we have not yet\ proved a time-independent positive lower bound of
the curvature.\ To prove it, we consider the case $0<\alpha<1\ $and $\alpha>1$ separately.

\subsubsection{The case $0<\alpha<1.$}

\ \ For the case $0<\alpha<1,$ we can improve Lemma \ref{lem-3} as:

\begin{lemma}
\label{lem-6}Under the LP\ flow (\ref{eqn1.3L}) on$\ S^{1}\times
\lbrack0,\infty)\ $with $0<\alpha<1,$ there exists a constant $c>0,\ $%
independent\ of time,\ such that%
\begin{equation}
k\left(  \theta,t\right)  \geq c>0,\ \ \ \forall\ \left(  \theta,t\right)  \in
S^{1}\times\lbrack0,\infty). \label{k-lower-bd-4}%
\end{equation}
$\ $The same result holds for the AP\ flow\ (\ref{eqn1.3A})\ on$\ S^{1}%
\times\lbrack0,\infty)\ $with $0<\alpha<1.$
\end{lemma}

%

\proof
We\ use the idea from p. 349 of\ \cite{And1}.\ In view of the relation
$k=d\theta/ds,\ $the gradient estimate in (\ref{kv}) can be rewritten as
$\left\vert \left(  k^{\alpha-1}\left(  \cdot,t\right)  \right)
_{s}\right\vert \leq C$ for some different\ constant $C\ $independent\ of
time,\ where$\ s\ $is the arc\ length parameter of $X\left(  \cdot,t\right)
.\ $Using this\ and the mean value theorem,\ together with the property that
the length of$\ X\left(  \cdot,t\right)  $ is decreasing in time in both
flows,\ we have
\begin{equation}
\left\vert \frac{1}{k^{1-\alpha}\left(  s_{2},t\right)  }-\frac{1}%
{k^{1-\alpha}\left(  s_{1},t\right)  }\right\vert \leq C\left\vert s_{2}%
-s_{1}\right\vert \leq C,\ \ \ 0<\alpha<1\label{kkcs}%
\end{equation}
for all $t\in\lbrack0,\infty)\ $and all $s_{1},\ s_{2}\ $on the curve
$X\left(  \cdot,t\right)  .\ $On the other hand, for each$\ t\in
\lbrack0,\infty),$ there exists some$\ \theta\left(  t\right)  \in S^{1}$ such
that
\[
L\left(  t\right)  =\int_{0}^{2\pi}\frac{1}{k\left(  \theta,t\right)  }%
d\theta=\frac{2\pi}{k\left(  \theta\left(  t\right)  ,t\right)  }\leq L\left(
0\right)  .
\]
This implies, at each time, the existence of some$\ $value of$\ s\left(
t\right)  $ such that
\begin{equation}
0<\frac{1}{k^{1-\alpha}\left(  s\left(  t\right)  ,t\right)  }\leq\left(
\frac{L\left(  0\right)  }{2\pi}\right)  ^{1-\alpha},\ \ \ 0<\alpha
<1.\label{kkcs-1}%
\end{equation}
The result follows by (\ref{kkcs-1})\ and\ (\ref{kkcs}).$%
\hfill
\square$

\ \ \ 

For the case $0<\alpha<1,\ $Theorem \ref{thm-1} and Lemma \ref{lem-6},
together with parabolic regularity theory, imply that both\ flows
are\ uniformly parabolic on $S^{1}\times\lbrack0,\infty)$\ and the evolving
curve $X\left(  \cdot,t\right)  $ converges\ as $t\rightarrow\infty$ in
$C^{\infty}\left(  S^{1}\right)  $ to a circle\ with radius $L\left(
0\right)  /2\pi\ $(in (\ref{eqn1.3L}))\ or $\sqrt{A\left(  0\right)  /\pi}%
\ $(in (\ref{eqn1.3A})). In fact, due to the uniform positive upper and lower
bounds of the curvature, Lemmas \ref{lem-7}, \ref{lem-8},\ \ref{lem-9}%
,\ Theorem \ref{thm-last} and Remark \ref{rmk1}\ in the next section\ are all
valid for the case $0<\alpha<1.$

Our conclusion is that now Theorem \ref{thm-main-L}\ and Theorem
\ref{thm-main-A} have been proved for the case $0<\alpha<1.\ $

\subsubsection{The case $\alpha>1.$}

It remains to prove a time-independent\ positive\ lower bound of the curvature
for the case $\alpha>1.\ $Note that the double integral in (\ref{1/k-2})\ may
tend to infinity as $t\rightarrow\infty.\ $Therefore, one needs to do more to
exclude the possibility that the lower bound may tend to zero as
$t\rightarrow\infty.\ $

We have:\ 

\begin{lemma}
\label{lem-7}Under the LP\ flow (\ref{eqn1.3L}) on$\ S^{1}\times
\lbrack0,\infty)\ $with $\alpha>1,\ $we have%
\begin{equation}
\frac{dA}{dt}\left(  t\right)  \rightarrow0\ \ \ \text{as \ \ }t\rightarrow
\infty\label{dAdt-0}%
\end{equation}
and under the AP\ flow (\ref{eqn1.3A}) on$\ S^{1}\times\lbrack0,\infty)\ $with
$\alpha>1,\ $we have%
\begin{equation}
\frac{dL}{dt}\left(  t\right)  \rightarrow0\ \ \ \text{as \ \ }t\rightarrow
\infty. \label{dLdt-0}%
\end{equation}

\end{lemma}%

\proof
Assume not for (\ref{dAdt-0}).\ Since $A^{\prime}\left(  t\right)  \geq
0,\ $there exists a constant $C_{0}>0\ $independent\ of time\ and a sequence
of times $\left\{  t_{i}\right\}  _{i=1}^{\infty}\ $going$\ $to infinity\ such
that $A^{\prime}\left(  t_{i}\right)  \geq C_{0}>0\ $for
all$\ i=1,\ 2,\ 3,\ ..$. Also note that%
\begin{align}
A^{\prime\prime}\left(  t\right)   &  =\frac{d}{dt}\left[  -\int_{0}^{2\pi
}k^{\alpha-1}\left(  \theta,t\right)  d\theta+\dfrac{{L\!}\left(  t\right)
}{2\pi}\int_{0}^{2\pi}k^{\alpha}\left(  \theta,t\right)  d\theta\right]
,\ \ \ L\left(  t\right)  =L\left(  0\right) \nonumber\\
&  =-\int_{0}^{2\pi}\left(  \alpha-1\right)  v\left(  v_{\theta\theta
}+v-\lambda\left(  t\right)  \right)  d\theta+\dfrac{{L\!}\left(  0\right)
}{2\pi}\int_{0}^{2\pi}\alpha v^{p}\left(  v_{\theta\theta}+v-\lambda\left(
t\right)  \right)  d\theta\nonumber\\
&  :=I\left(  t\right)  +II\left(  t\right)  , \label{second-A}%
\end{align}
where $p=1+1/\alpha.\ $By (\ref{kv}),\ we have%
\[
\left\vert I\left(  t\right)  \right\vert =\left\vert -\left(  \alpha
-1\right)  \int_{0}^{2\pi}\left(  -v_{\theta}^{2}+v^{2}-v\lambda\left(
t\right)  \right)  d\theta\right\vert \leq C_{1},\ \ \ \lambda\left(
t\right)  =\dfrac{1}{{2\pi}}%
{\displaystyle\int_{0}^{2\pi}}
k^{\alpha}\left(  \theta,t\right)  d\theta
\]
for some constant $C_{1}\ $independent\ of time. Also
\begin{equation}
\left\vert II\left(  t\right)  \right\vert =\left\vert \dfrac{{L\!}\left(
0\right)  }{2\pi}\int_{0}^{2\pi}\left(  -\alpha pv^{p-1}v_{\theta}^{2}+\alpha
v^{p+1}-\alpha v^{p}\lambda\left(  t\right)  \right)  d\theta\right\vert \leq
C_{2},\ \ \ p=1+\frac{1}{\alpha}>1 \label{2000}%
\end{equation}
for another constant $C_{2}\ $independent\ of time. Therefore we obtain
$\left\vert A^{\prime\prime}\left(  t\right)  \right\vert \leq C_{3}\ $for all
$t\in\lbrack0,\infty)$ for some constant $C_{3}\ $independent\ of time. As the
derivative of $A^{\prime}\left(  t\right)  \ $is uniformly bounded, one can
find a number $\rho_{0}$ independent\ of $t_{i}$ such that%
\[
A^{\prime}\left(  t\right)  \geq\frac{C_{0}}{2}>0,\ \ \ \forall\ t\in\left[
t_{i},t_{i}+\rho_{0}\right]  ,\ \ \ i=1,\ 2,\ 3,\ ....
\]
This implies that $A\left(  \infty\right)  -A\left(  0\right)  =\int
_{0}^{\infty}A^{\prime}\left(  t\right)  dt=\infty,$ contradicting to the
inequality $0<4\pi A\left(  t\right)  \leq L^{2}\left(  t\right)
=L^{2}\left(  0\right)  \ $for all time.\ The proof of\ (\ref{dAdt-0})\ is done.

Similarly, assume not for (\ref{dLdt-0}).\ Since $L^{\prime}\left(  t\right)
\leq0,\ $there exists a constant $C_{0}>0\ $and a sequence of times $\left\{
t_{i}\right\}  _{i=1}^{\infty}\ $going$\ $to infinity\ such that $L^{\prime
}\left(  t_{i}\right)  \leq-C_{0}<0$ for all$\ i=1,\ 2,\ 3,\ ..$. Also note
that
\begin{align}
L^{\prime\prime}\left(  t\right)   &  =\frac{d}{dt}\left[  -\int_{0}^{2\pi
}k^{\alpha}d\theta+\dfrac{2\pi}{{L}\left(  t\right)  }\int_{0}^{2\pi}%
k^{\alpha-1}d\theta\right] \nonumber\\
&  =\left\{
\begin{array}
[c]{l}%
-%
{\displaystyle\int_{0}^{2\pi}}
\alpha v^{p}\left(  v_{\theta\theta}+v-\lambda\left(  t\right)  \right)
d\theta+\dfrac{2\pi}{{L}\left(  t\right)  }%
{\displaystyle\int_{0}^{2\pi}}
\left(  \alpha-1\right)  v\left(  v_{\theta\theta}+v-\lambda\left(  t\right)
\right)  d\theta%
\vspace{3mm}%
\\
-\dfrac{2\pi}{L^{2}\left(  t\right)  }\left(  -%
{\displaystyle\int_{0}^{2\pi}}
k^{\alpha}d\theta+\dfrac{2\pi}{{L}\left(  t\right)  }%
{\displaystyle\int_{0}^{2\pi}}
k^{\alpha-1}d\theta\right)
{\displaystyle\int_{0}^{2\pi}}
k^{\alpha-1}d\theta,\ \ \ \alpha>1,
\end{array}
\right.  \label{second-L}%
\end{align}
where $\lambda\left(  t\right)  =$ ${L}^{-1}\left(  t\right)  \int_{0}^{2\pi
}k^{\alpha-1}\left(  \theta,t\right)  d\theta,$ $\alpha>1.$ By (\ref{kv}%
)\ again, $\lambda\left(  t\right)  $ is uniformly bounded and we have
$\left\vert L^{\prime\prime}\left(  t\right)  \right\vert \leq C_{4}$ for all
$t\in\lbrack0,\infty)$ for some constant $C_{4}\ $independent\ of time. As the
derivative of $L^{\prime}\left(  t\right)  \ $is uniformly bounded, one can
find a number $\rho_{0}$ independent\ of $t_{i}$ such that%
\[
L^{\prime}\left(  t\right)  \leq-\frac{C_{0}}{2}<0,\ \ \ \forall\ t\in\left[
t_{i},t_{i}+\rho_{0}\right]  ,\ \ \ i=1,\ 2,\ 3,\ ....
\]
This implies that $L\left(  \infty\right)  -L\left(  0\right)  =\int
_{0}^{\infty}L^{\prime}\left(  t\right)  dt=-\infty,$ contradicting to the
isoperimetric inequality $L^{2}\left(  t\right)  \geq4\pi A\left(  t\right)
=4\pi A\left(  0\right)  \ $for all time.\ The proof is done.$%
\hfill
\square\ $\ 

\ \ \ 

To go further we recall the following inequality\ in\ Andrews\ \cite{And1}:

\begin{lemma}
Let $M$ be a compact Riemannian manifold with a volume form $d\mu,$\ and let
$\xi$\ be a continuous function on$\;M.$\ Then for any\textbf{ }%
\emph{decreasing} continuous function $F:\mathbb{R\rightarrow R},\;$we have
\begin{equation}
\int_{M}\xi d\mu\int_{M}F\left(  \xi\right)  d\mu\geq\int_{M}d\mu\int_{M}\xi
F\left(  \xi\right)  d\mu. \label{Ben-Andrews}%
\end{equation}
If $F$\ is strictly decreasing, then equality holds if and only if $\xi$\ is a
constant function on $M.$\ 
\end{lemma}

\begin{remark}
If $F$ is increasing, then we have $\leq$ in (\ref{Ben-Andrews}).
\end{remark}

Next, we prove the following:\ 

\begin{lemma}
\label{lem-8}Under the LP\ flow\ (\ref{eqn1.3L}) on$\ S^{1}\times
\lbrack0,\infty)\ $with $\alpha>1,\ $we have%
\begin{equation}
\lim_{t\rightarrow\infty}\left\Vert k\left(  \cdot,t\right)  -\frac{2\pi
}{L\left(  0\right)  }\right\Vert _{C^{0}\left(  S^{1}\right)  }=0.
\label{conv1}%
\end{equation}

\end{lemma}%

\proof
Assume not.\ Then there exists a sequence of times $\left\{  t_{i}\right\}
_{i=1}^{\infty}\ $going$\ $to infinity\ such that
\begin{equation}
\left\Vert k\left(  \cdot,t_{i}\right)  -\frac{2\pi}{L\left(  0\right)
}\right\Vert _{C^{0}\left(  S^{1}\right)  }\geq\varepsilon>0,\ \ \ \forall\ i
\label{con}%
\end{equation}
for some $\varepsilon>0.\ $As $v$ and $v_{\theta}$ are both uniformly
bounded,\ there is a subsequence of $\left\{  t_{i}\right\}  _{i=1}^{\infty
},\ $still denoted as $\left\{  t_{i}\right\}  _{i=1}^{\infty},\ $such
that$\ v\left(  \theta,t_{i}\right)  =k^{\alpha}\left(  \theta,t_{i}\right)
\ $converges uniformly on\ $S^{1}$ to a Lipschitz continuous function
$w(\theta)\geq0$ as $i\rightarrow\infty.\ $In particular,\ since $\alpha
>1,\ $both $k\left(  \theta,t_{i}\right)  \ $and$\ k^{\alpha-1}\left(
\theta,t_{i}\right)  $\ converge uniformly to $w^{1/\alpha}\left(
\theta\right)  \ $and $w^{1-1/\alpha}\left(  \theta\right)  $%
\ respectively.$\ $At\ this moment, we can not guarantee that $w^{1/\alpha
}\left(  \theta\right)  \geq0\ $is strictly positive.\ Hence we cannot
conclude that the integral of$\ 1/k\left(  \theta,t_{i}\right)  \ $converges
to the integral of$\ w^{-1/\alpha}\left(  \theta\right)  .\ $One needs to do
more. By (\ref{dAdt-0}), we have$\ $%
\begin{equation}
0=\lim_{i\rightarrow\infty}\frac{dA}{dt}\left(  t_{i}\right)  =-\int_{0}%
^{2\pi}w^{1-\frac{1}{\alpha}}\left(  \theta\right)  d\theta+\dfrac
{{L\!}\left(  0\right)  }{2\pi}\int_{0}^{2\pi}w\left(  \theta\right)
d\theta,\ \ \ 1-\frac{1}{\alpha}>0. \label{dAdt}%
\end{equation}
Also note that
\[
L\left(  0\right)  =\int_{0}^{2\pi}\frac{1}{k\left(  \theta,t_{i}\right)
}d\theta,\ \ \ \forall\ i
\]
and by Fatou's lemma, we have%
\begin{equation}
0\leq\int_{0}^{2\pi}\frac{1}{w^{1/\alpha}\left(  \theta\right)  }d\theta\leq
L\left(  0\right)  . \label{fatou}%
\end{equation}
Thus the Lebesgue\ integral $\int_{0}^{2\pi}w^{-1/\alpha}\left(
\theta\right)  d\theta\ $converges,\ i.e.,\ $w^{1/\alpha}\left(
\theta\right)  >0\ $almost everywhere\ on\ $S^{1}$.\ Now by (\ref{dAdt})\ and
(\ref{fatou}), we have
\begin{equation}
2\pi\int_{0}^{2\pi}w^{1-\frac{1}{\alpha}}\left(  \theta\right)  d\theta
\geq\int_{0}^{2\pi}\frac{1}{w^{1/\alpha}\left(  \theta\right)  }d\theta
\int_{0}^{2\pi}w\left(  \theta\right)  d\theta. \label{ineq1}%
\end{equation}
On the other hand, by inequality (\ref{Ben-Andrews}), we also have (take
$F\left(  \xi\right)  =1/\left(  \xi^{1/\alpha}+\varepsilon\right)  ,\ \xi
\geq0,\ \varepsilon>0$)%
\[
\int_{0}^{2\pi}w\left(  \theta\right)  d\theta\int_{0}^{2\pi}\frac
{1}{w^{1/\alpha}\left(  \theta\right)  +\varepsilon}d\theta\geq2\pi\int
_{0}^{2\pi}\frac{w\left(  \theta\right)  }{w^{1/\alpha}\left(  \theta\right)
+\varepsilon}d\theta,\ \ \ \forall\ \varepsilon>0.
\]
Letting the constant $\varepsilon\rightarrow0^{+},\ $we obtain%
\begin{equation}
\int_{0}^{2\pi}w\left(  \theta\right)  d\theta\int_{0}^{2\pi}\frac
{1}{w^{1/\alpha}\left(  \theta\right)  }d\theta\geq2\pi\int_{0}^{2\pi
}w^{1-\frac{1}{\alpha}}\left(  \theta\right)  d\theta. \label{ineq2}%
\end{equation}
Combine (\ref{ineq1})\ and (\ref{ineq2}) to get%
\[
2\pi\int_{0}^{2\pi}w^{1-\frac{1}{\alpha}}\left(  \theta\right)  d\theta
=\int_{0}^{2\pi}\frac{1}{w^{1/\alpha}\left(  \theta\right)  }d\theta\int
_{0}^{2\pi}w\left(  \theta\right)  d\theta,
\]
i.e., the following iterated integral vanishes
\begin{equation}
\frac{1}{2}\int_{0}^{2\pi}\int_{0}^{2\pi}\left[  \left(  \frac{1}{w^{1/\alpha
}\left(  x\right)  }-\frac{1}{w^{1/\alpha}\left(  y\right)  }\right)  \left(
w\left(  x\right)  -w\left(  y\right)  \right)  \right]  dxdy=0. \label{ineq3}%
\end{equation}
As the integrand in (\ref{ineq3})\ is nonpositive almost
everywhere\ on\ $S^{1}\times S^{1}\ $and$\ w\left(  \theta\right)  \geq0\ $is
a continuous function,$\ $we must have$\ w\left(  \theta\right)
=C\ $everywhere on $S^{1},$ where $C\ $is the constant given by$\ \left(
2\pi/L\left(  0\right)  \right)  ^{\alpha}\ $due\ to (\ref{dAdt}). This
implies that$\ k\left(  \cdot,t_{i}\right)  \ $converges uniformly to the
constant$\ 2\pi/L\left(  0\right)  \ $as $i\rightarrow\infty,\ $contradicting
to (\ref{con}).\ The proof is done.$%
\hfill
\square$

\ 

Finally, we have:

\begin{lemma}
\label{lem-9}Under the AP\ flow\ (\ref{eqn1.3A}) on$\ S^{1}\times
\lbrack0,\infty)\ $with $\alpha>1,\ $we have%
\begin{equation}
\lim_{t\rightarrow\infty}\left\Vert k\left(  \cdot,t\right)  -\sqrt{\frac{\pi
}{A\left(  0\right)  }}\right\Vert _{C^{0}\left(  S^{1}\right)  }=0.
\label{conv2}%
\end{equation}

\end{lemma}%

\proof
The proof is also identical to Lemma \ref{lem-8}.\ We now have$\ $
\begin{equation}
0=\lim_{i\rightarrow\infty}\frac{dL}{dt}\left(  t_{i}\right)  =-\int_{0}%
^{2\pi}w\left(  \theta\right)  d\theta+\dfrac{2\pi}{{L\!}\left(
\infty\right)  }\int_{0}^{2\pi}w^{1-\frac{1}{\alpha}}\left(  \theta\right)
d\theta,\ \ \ 1-\frac{1}{\alpha}>0. \label{dLdt}%
\end{equation}
where ${L}\left(  \infty\right)  =\lim_{t\rightarrow\infty}L\left(  t\right)
\ $and\ by Fatou's lemma, we have (\ref{fatou}) and so we arrive at the
inequality (\ref{ineq1}) again. The rest of the proof goes the same as in
Lemma \ref{lem-8} and we conclude that $k\left(  \cdot,t\right)  \ $converges
uniformly to the constant$\ 2\pi/L\left(  \infty\right)  \ $as $t\rightarrow
\infty.\ $Since\ the limit is a circle, we have$\ 2\pi/L\left(  \infty\right)
=\sqrt{\pi/A\left(  \infty\right)  }=\sqrt{\pi/A\left(  0\right)  }.\ $The
proof is done.$%
\hfill
\square$

\ \ \ 

Now by Lemma \ref{lem-8} and Lemma \ref{lem-9}\ we can obtain the following
$C^{\infty}$ convergence of curvature:

\begin{theorem}
\label{thm-last}Under the\ LP\ flow\ (\ref{eqn1.3L}) with $\alpha>1\ $we have%
\begin{equation}
\lim_{t\rightarrow\infty}\left\Vert k\left(  \cdot,t\right)  -\frac{2\pi
}{L\left(  0\right)  }\right\Vert _{C^{m}\left(  S^{1}\right)  }%
=0,\ \ \ \forall\ m=0,\ 1,\ 2,\ 3,\ .... \label{final-1}%
\end{equation}
and under the AP flow\ (\ref{eqn1.3A}) with $\alpha>1$ we have%
\begin{equation}
\lim_{t\rightarrow\infty}\left\Vert k\left(  \cdot,t\right)  -\sqrt{\frac{\pi
}{A\left(  0\right)  }}\right\Vert _{C^{m}\left(  S^{1}\right)  }%
=0,\ \ \ \forall\ m=0,\ 1,\ 2,\ 3,\ .... \label{final-2}%
\end{equation}

\end{theorem}%

\proof
By (\ref{conv1}), the curvature equation in (\ref{dkdt}) is uniformly
parabolic and the parabolic regularity theory implies that all space-time
derivatives of$\ k\left(  \theta,t\right)  \ $are uniformly\ bounded by
constants depending only on the order of differentiation. By induction and the
Arzela-Ascoli theorem, we can easily obtain (\ref{final-1}) and (\ref{final-2}%
).$\
\hfill
\square$

\ \ 

The proof of Theorem \ref{thm-main-L} and\ Theorem \ref{thm-main-A}\ is now complete.\ 

\begin{remark}
\label{rmk1}Strictly speaking, the proof of Theorem \ref{thm-main-L}
and\ Theorem \ref{thm-main-A}\ is not quite complete yet because in Theorem
\ref{thm-last} it is only shown that, in both flows, the curvature tends to a
constant eventually. Conceivably the evolving curve $X\left(  \cdot,t\right)
\ $may escape to infinity or oscillate indefinitely\footnote{We thank the
referee for pointing out this important issue, which we had neglected at
first.\ \ }. To show that this will not happen, one needs additional arguments
to ensure that the flow does tend to a unique fixed\ circle eventually. We
only give a brief description here.\ The idea is to prove that the curvature
$k\left(  \cdot,t\right)  $ converges to a constant in an
\emph{exponentially-decay} way, like what has been done in
Gage-Hamilton\ \cite{GH}.\ One can also see Theorem 5.10. in the recent paper
\cite{MPW}, which has successfully generalized Gage-Hamilton's argument to the
nonlocal flow (\ref{eqn1.3A}) for the case $\alpha\in\mathbb{N}.\ $Recall that
the function$\ v=k^{\alpha}$ satisfies the equation%
\begin{equation}
v_{t}\left(  \theta,t\right)  =\alpha v^{p}\left(  \theta,t\right)  \left[
v_{\theta\theta}\left(  \theta,t\right)  +v\left(  \theta,t\right)
-\lambda\left(  t\right)  \right]  ,\ \ \ p=1+\frac{1}{\alpha}>1,\ \ \ \left(
\theta,t\right)  \in S^{1}\times\lbrack0,\infty) \label{v-eq-2}%
\end{equation}
and we have know that $v\left(  \theta,t\right)  $ tends to a
positive\ constant as $t\rightarrow\infty\ $in the space$\ C^{\infty}\left(
S^{1}\right)  $ and $\lambda\left(  t\right)  $ also tends to a
positive\ constant as $t\rightarrow\infty$ (with all of its derivatives tend
to zero\ as $t\rightarrow\infty$). With this, plus a further analysis (for
example,\ \emph{linearization}) of the equation (\ref{v-eq-2}), can convince
us that the difference $v\left(  \theta,t\right)  -\lambda\left(  t\right)  $
does decay to zero exponentially\ in $C^{\infty}\left(  S^{1}\right)  $\ (one
can also mimic the method in \cite{MPW} or see \cite{ES}).$\ $Next, by the
equation of the support function, which is\ (see (\ref{supp}))%
\begin{equation}
u_{t}\left(  \theta,t\right)  =-v\left(  \theta,t\right)  +\lambda\left(
t\right)  =-k^{\alpha}\left(  \theta,t\right)  +\lambda\left(  t\right)
,\ \ \ \left(  \theta,t\right)  \in S^{1}\times\lbrack0,\infty),
\end{equation}
we know $u_{t}\left(  \theta,t\right)  \ $also\ decays to zero exponentially
in $C^{\infty}\left(  S^{1}\right)  .$ By the identity$\ u_{\theta\theta
}+u=1/k,\ $one can express$\ u$ as
\begin{equation}
u\left(  \theta,t\right)  =a\left(  t\right)  \cos\theta+b\left(  t\right)
\sin\theta+\int_{0}^{\theta}\frac{1}{k\left(  \sigma,t\right)  }\sin\left(
\theta-\sigma\right)  d\sigma,\ \ \left(  \theta,t\right)  \in S^{1}%
\times\lbrack0,\infty), \label{u-exp}%
\end{equation}
where $a\left(  t\right)  =u\left(  0,t\right)  \ $and$\ b\left(  t\right)
=u_{\theta}\left(  0,t\right)  .$ By (\ref{u-exp}), we get%
\[
u_{\theta}\left(  \theta,t\right)  =-a\left(  t\right)  \sin\theta+b\left(
t\right)  \cos\theta+\int_{0}^{\theta}\frac{1}{k\left(  \sigma,t\right)  }%
\cos\left(  \theta-\sigma\right)  d\sigma,
\]
and then%
\[
u_{t}\left(  0,t\right)  =a^{\prime}\left(  t\right)  ,\ \ \ \left(
u_{\theta}\right)  _{t}\left(  0,t\right)  =\left(  u_{t}\right)  _{\theta
}\left(  0,t\right)  =b^{\prime}\left(  t\right)  ,\ \ \ t\in\lbrack
0,\infty).
\]
Since $a^{\prime}\left(  t\right)  \ $and $b^{\prime}\left(  t\right)  $ both
decay to zero exponentially, the two integrals $\int_{0}^{\infty}a^{\prime
}\left(  t\right)  dt\ $and $\int_{0}^{\infty}b^{\prime}\left(  t\right)  dt$
both converge.\ By%
\[
\lim_{t\rightarrow\infty}a\left(  t\right)  =a\left(  0\right)  +\int
_{0}^{\infty}a^{\prime}\left(  t\right)  dt,\ \ \ \lim_{t\rightarrow\infty
}b\left(  t\right)  =b\left(  0\right)  +\int_{0}^{\infty}b^{\prime}\left(
t\right)  dt
\]
and (\ref{u-exp}),\ we conclude that the support function$\ u$ has a limit:
\begin{equation}
\lim_{t\rightarrow\infty}u\left(  \theta,t\right)  =a\cos\theta+b\sin
\theta+c,\ \ \ \forall\ \theta\in S^{1}, \label{u-ABC}%
\end{equation}
for some constants $a,\ b,\ c\in\mathbb{R\ }$with $c>0.\ $Thus the circle in
Theorem \ref{thm-main-L} is a fixed circle centered at $\left(  a,b\right)
\in\mathbb{R}^{2}$ with radius $L\left(  0\right)  /2\pi\ $and the circle in
Theorem \ref{thm-main-A} is a fixed circle centered at $\left(  a,b\right)
\in\mathbb{R}^{2}$ with radius $\sqrt{A\left(  0\right)  /\pi}.$
\end{remark}

\section{Other $k^{\alpha}$-type nonlocal flow for $\alpha\geq1$.
\label{general}}

When$\ \alpha\geq1,\ $using certain geometric inequalities for convex closed
plane curves, we can consider two additional types of nonlocal flow.\ Both
flows have the interesting feature that $L\left(  t\right)  \ $is decreasing
and $A\left(  t\right)  \ $is\ increasing.$\ $

Recall the following two inequalities for convex closed plane curves\ (note
that (\ref{11})\ is the same as (\ref{kkL})):%
\begin{equation}
\int_{0}^{2\pi}k^{\beta-1}\left(  \theta\right)  d\theta\leq\frac{L}{2\pi}%
\int_{0}^{2\pi}k^{\beta}\left(  \theta\right)  d\theta,\ \ \ \forall
\text{\ }\beta\geq0 \label{11}%
\end{equation}
and
\begin{equation}
\int_{0}^{2\pi}k^{\beta}\left(  \theta\right)  d\theta\leq\frac{2A}{L}\int
_{0}^{2\pi}k^{\beta+1}\left(  \theta\right)  d\theta,\ \ \ \forall
\text{\ }\beta\geq0, \label{22}%
\end{equation}
where for $\beta>0\ $the equality holds in\ (\ref{11})\ if and only if
$\gamma$ is\ a circle; and for $\beta\geq0\ $the equality holds in\ (\ref{22}%
)\ if and only if $\gamma$ is\ a circle.\ When $\beta=0,\ $we have equality
in\ (\ref{11}) for any convex closed curve; and when $\beta=0,\ $(\ref{22})
becomes%
\begin{equation}
\int_{0}^{2\pi}k\left(  \theta\right)  d\theta\geq\frac{\pi L}{A},
\label{gage}%
\end{equation}
which is precisely\ \emph{Gage's isoperimetric\ inequality} for convex closed
curves (see \cite{Ga3}).\ In fact, (\ref{22})\ is a combination of Gage's
inequality\ and\ H\"{o}lder inequality\ (or Andrews's inequality
(\ref{Ben-Andrews})).\ One can also see Corollary 1.10. in \cite{LT}.

Combined with (\ref{22})\ and the inequality $L^{2}\geq4\pi A,$ for
$\alpha\geq1\ $we have the following two inequalities for any convex closed
plane curve:%
\begin{equation}
\frac{1}{L}\int_{0}^{2\pi}k^{\alpha-1}d\theta\leq\frac{2A}{L^{2}}\int
_{0}^{2\pi}k^{\alpha}d\theta\leq\frac{1}{2\pi}\int_{0}^{2\pi}k^{\alpha}%
d\theta,\ \ \ \alpha\geq1 \label{gage1}%
\end{equation}
or%
\begin{equation}
\frac{1}{L}\int_{0}^{2\pi}k^{\alpha-1}d\theta\leq\frac{L}{4\pi A}\int
_{0}^{2\pi}k^{\alpha-1}d\theta\leq\frac{1}{2\pi}\int_{0}^{2\pi}k^{\alpha
}d\theta,\ \ \ \alpha\geq1. \label{gage2}%
\end{equation}
Motivated by (\ref{gage1})\ and (\ref{gage2}), we can prove the following:

\begin{theorem}
\label{thm-2}Assume$\ \alpha\geq1\ $and $X_{0}\left(  u\right)  ,\ u\in
S^{1},$ is a smooth convex closed curve. Then the nonlocal flow\ (\ref{eqn1.2}%
), where
\begin{equation}
F\left(  k\right)  =k^{\alpha},\ \ \ \lambda\left(  t\right)  =\frac{2A\left(
t\right)  }{L^{2}\left(  t\right)  }\int_{X\left(  \cdot,t\right)  }%
k^{\alpha+1}ds, \label{g1}%
\end{equation}
has a smooth convex solution for all time $t\in\lbrack0,\infty).$ Moreover, it
converges to a round circle in $C^{\infty}\ $topology\ as $t\rightarrow
\infty.$ The same result holds if $\ $%
\begin{equation}
F\left(  k\right)  =k^{\alpha},\ \ \ \lambda\left(  t\right)  =\frac{L\left(
t\right)  }{4\pi A\left(  t\right)  }\int_{X\left(  \cdot,t\right)  }%
k^{\alpha}ds. \label{g2}%
\end{equation}

\end{theorem}

\begin{remark}
When $\alpha=1,\ \lambda\left(  t\right)  $ in (\ref{g2})\ becomes$\ \lambda
=L/2A,\ $which is a\ gradient flow of the isoperimetric ratio functional.\ It
has been studied in \cite{JP}.
\end{remark}

\begin{remark}
By (\ref{gage1}), one can view (\ref{g1}) as a flow below\ the LP case (hence
$L\left(  t\right)  \ $is decreasing) and above the the AP\ case (hence
$A\left(  t\right)  \ $is increasing).\ The same for (\ref{g2}).
\end{remark}

%

\proof
Again, the flow has a unique smooth convex solution for short time\ $[0,T)\ $%
and we can use the outward normal angle $\theta$\ as\ a parameter.\ Now for
the flow (\ref{eqn1.2})\ with $F\left(  k\right)  \ $and$\ \lambda\left(
t\right)  \ $given by (\ref{g1})\ or (\ref{g2}), we have \
\begin{equation}
\frac{dL}{dt}=2\pi\lambda\left(  t\right)  -\int_{0}^{2\pi}k^{\alpha}\left(
\theta,t\right)  d\theta\leq0,\ \ \ t\in\lbrack0,T)
\end{equation}
and%
\begin{equation}
\frac{dA}{dt}=\lambda\left(  t\right)  L\left(  t\right)  -\int_{0}^{2\pi
}k^{\alpha-1}\left(  \theta,t\right)  d\theta\geq0,\ \ \ t\in\lbrack0,T).
\end{equation}
Thus $L\left(  t\right)  \ $is decreasing and $A\left(  t\right)
\ $is\ increasing.\ In particular, the isoperimetric ratio $L^{2}\left(
t\right)  /4\pi A\left(  t\right)  $ is decreasing.\ 

Next, we claim that all of the previous curvature estimates are valid for
$\lambda\left(  t\right)  $ given by (\ref{g1})\ or (\ref{g2}). First, the
proof in Lemma \ref{lem-1} (i.e., Tso's estimate) carries over to this case
since now the length is decreasing and the area is increasing\ (in both
(\ref{g1}) and (\ref{g2})).\ Moreover, the inequality (\ref{Qt})\ is valid as
long as $\lambda\left(  t\right)  \ $is a positive
quantity.\ Similarly,\ Lemma \ref{lem-2} is also valid as long as
$\lambda\left(  t\right)  \ $is positive (see Remark \ref{rmk1}).\ In the
proof of inequality (\ref{phi}), the term $\lambda\left(  t\right)  $ has been
cancelled.\ Hence Lemma \ref{lem-3} is also valid here.\ Consequently, we have
Theorem \ref{thm-1} and the flow is smooth, convex, defined on $S^{1}%
\times\lbrack0,\infty).$

Next we check the validity of Lemma \ref{lem-7}.\ For $\lambda\left(
t\right)  $ given by (\ref{g1}) we have
\[
\lambda^{\prime}\left(  t\right)  =\frac{2A\left(  t\right)  }{L^{2}\left(
t\right)  }\int_{0}^{2\pi}\alpha v^{p}\left(  v_{\theta\theta}+v-\lambda
\left(  t\right)  \right)  d\theta+\left(  \int_{0}^{2\pi}k^{\alpha}%
d\theta\right)  \frac{d}{dt}\left(  \frac{2A\left(  t\right)  }{L^{2}\left(
t\right)  }\right)  ,\ \ \ p=1+\frac{1}{\alpha},
\]
and, same as before, it is easy to see that $\left\vert \lambda\left(
t\right)  \right\vert \leq C\ $and$\ \left\vert \lambda^{\prime}\left(
t\right)  \right\vert \leq C$ for all $t\in\lbrack0,\infty)$ for some constant
$C$ independent\ of time. Hence we have \
\[
\left\vert \frac{d^{2}A}{dt^{2}}\left(  t\right)  \right\vert =\left\vert
\lambda^{\prime}\left(  t\right)  L\left(  t\right)  +\lambda\left(  t\right)
L^{\prime}\left(  t\right)  -\int_{0}^{2\pi}\left(  \alpha-1\right)  v\left(
v_{\theta\theta}+v-\lambda\left(  t\right)  \right)  d\theta\right\vert \leq
C
\]
and
\[
\left\vert \frac{d^{2}L}{dt^{2}}\left(  t\right)  \right\vert =\left\vert
2\pi\lambda^{\prime}\left(  t\right)  -\int_{0}^{2\pi}\alpha v^{p}\left(
v_{\theta\theta}+v-\lambda\left(  t\right)  \right)  d\theta\right\vert \leq
C
\]
for all $t\in\lbrack0,\infty)$ for some constant $C$ independent\ of time. The
above two estimates imply
\begin{equation}
\frac{dA}{dt}\left(  t\right)  =\lambda\left(  t\right)  L\left(  t\right)
-\int_{0}^{2\pi}k^{\alpha-1}\left(  \theta,t\right)  d\theta\rightarrow
0\ \ \ \text{as\ \ \ }t\rightarrow\infty\label{dAdt-000}%
\end{equation}
and%
\begin{equation}
\frac{dL}{dt}\left(  t\right)  =2\pi\lambda\left(  t\right)  -\int_{0}^{2\pi
}k^{\alpha}\left(  \theta,t\right)  d\theta\rightarrow0\ \ \ \text{as\ \ \ }%
t\rightarrow\infty. \label{dLdt-000}%
\end{equation}
Thus Lemma \ref{lem-7} is valid for $\lambda\left(  t\right)  $ given by
(\ref{g1}).\ One can check that Lemma \ref{lem-7} is also valid for\ $\lambda
\left(  t\right)  $ given by (\ref{g2}),

Finally, it remains to check the$\ C^{0}$ convergence of the curvature.\ For
$\lambda\left(  t\right)  $ given by (\ref{g1}), since we have both
(\ref{dAdt-000})\ and (\ref{dLdt-000}),\ (\ref{dAdt})\ in Lemma\ \ref{lem-8}
can be replaced by
\begin{equation}
0=\lim_{i\rightarrow\infty}\frac{dL}{dt}\left(  t_{i}\right)  =\frac{4\pi
A\left(  \infty\right)  }{L^{2}\left(  \infty\right)  }\int_{0}^{2\pi}w\left(
\theta\right)  d\theta-\int_{0}^{2\pi}w\left(  \theta\right)  d\theta
\end{equation}
and%
\begin{equation}
0=\lim_{i\rightarrow\infty}\frac{dA}{dt}\left(  t_{i}\right)  =\frac{2A\left(
\infty\right)  }{L\left(  \infty\right)  }\int_{0}^{2\pi}w\left(
\theta\right)  d\theta-\int_{0}^{2\pi}w^{1-\frac{1}{\alpha}}\left(
\theta\right)  d\theta. \label{s1}%
\end{equation}
Hence$\ L^{2}\left(  \infty\right)  =4\pi A\left(  \infty\right)  $ and
(\ref{s1}) becomes
\[
\frac{L\left(  \infty\right)  }{2\pi}\int_{0}^{2\pi}w\left(  \theta\right)
d\theta=\int_{0}^{2\pi}w^{1-\frac{1}{\alpha}}\left(  \theta\right)  d\theta.
\]
By Fatou's lemma again, we get%
\[
0\leq\int_{0}^{2\pi}\frac{1}{w^{1/\alpha}\left(  \theta\right)  }d\theta\leq
L\left(  \infty\right)
\]
and conclude the same inequality as\ (\ref{ineq1}). Same as in the proof of
Lemma \ref{lem-8}$\ $we have$\ w\left(  \theta\right)  =C\ $everywhere on
$S^{1}\ $for some constant $C,\ $where $C>0\ $is the constant given
by$\ \left(  2\pi/L\left(  \infty\right)  \right)  ^{\alpha}.\ $This implies
that$\ k\left(  \cdot,t_{i}\right)  \ $converges uniformly to the
constant$\ 2\pi/L\left(  \infty\right)  \ $as $i\rightarrow\infty.$ The proof
for the case (\ref{g2}) is similar.\ 

Until now, we have obtained a $C^{0}$ convergence of the curvature of the flow
for $\lambda\left(  t\right)  $ given by (\ref{g1})\ or (\ref{g2}).\ The
$C^{\infty}$ convergence of the flow follows from standard regularity
theory.\ The proof of Theorem \ref{thm-2} is finished (see Remark \ref{rmk1}
also).$%
\hfill
\square$

\section{Appendix; entropy estimate.}

Results in this short appendix section are supplementary.\ They are
interesting on their own. The property that we will claim is the so-called
"\emph{entropy estimate}" (we take this terminology from \cite{And1},\ p.
322). It is about the monotonicity of certain integral quantities.

We need the following inequality,\ which is a consequence of the Minkowski
inequality (1.6)\ in p. 322 of\ \cite{And1}:

\begin{lemma}
\label{lem-Mink}Let $\gamma$ be a convex closed curve with curvature $k\left(
\theta\right)  >0,\ $where $\theta$ is its outward normal angle. Then\ for any
$C^{2}\ $function $\Phi:\left(  0,\infty\right)  \rightarrow\mathbb{R},$ we
have the following inequalities:%
\begin{equation}
2\pi\int_{0}^{2\pi}\Phi\left(  k\right)  \left[  \Phi\left(  k\right)
_{\theta\theta}+\Phi\left(  k\right)  \right]  d\theta\leq\left(  \int
_{0}^{2\pi}\Phi\left(  k\right)  d\theta\right)  ^{2} \label{in-1}%
\end{equation}
and%
\begin{equation}
2A\int_{0}^{2\pi}\Phi\left(  k\right)  \left[  \Phi\left(  k\right)
_{\theta\theta}+\Phi\left(  k\right)  \right]  d\theta\leq\left(  \int
_{0}^{2\pi}\frac{\Phi\left(  k\right)  }{k}d\theta\right)  ^{2}, \label{in-2}%
\end{equation}
where $\Phi\left(  k\right)  =\Phi\left(  k\left(  \theta\right)  \right)
\ $and$\ A$ is the enclosed area of $\gamma.$
\end{lemma}

%

\proof
In\ p. 322 of \cite{And1}, for a convex closed plane\ curve $\gamma\ $with
support function $u\left(  \theta\right)  \ $and curvature $k\left(
\theta\right)  $, we can take $s_{1}\left(  \theta\right)  =u\left(
\theta\right)  \ $and $s_{2}\left(  \theta\right)  =\Phi\left(  k\left(
\theta\right)  \right)  \ $in the Minkowski inequality,\ where$\ \Phi\left(
z\right)  \ $is an arbitrary $C^{2}\ $function on\ $(0,\infty)$.\ This proves
(\ref{in-2}).\ On the other hand, one can also take $s_{1}\left(
\theta\right)  =1\ $and $s_{2}\left(  \theta\right)  =\Phi\left(  k\left(
\theta\right)  \right)  $. Then we get\ (\ref{in-1}).\ $%
\hfill
\square$

\ 

If we consider the flow$\ {X}_{t}=\Phi\left(  k\right)  \mathbf{N}_{in}{,\ }%
$then the length $L\left(  t\right)  $ and enclosed area $A\left(  t\right)  $
of $X\left(  \cdot,t\right)  $ will satisfy%
\begin{equation}
\frac{dA}{dt}\left(  t\right)  =-\int_{0}^{2\pi}\frac{\Phi\left(  k\left(
\theta,t\right)  \right)  }{k\left(  \theta,t\right)  }d\theta,\ \ \ \ \ \frac
{dL}{dt}\left(  t\right)  =-\int_{0}^{2\pi}\Phi\left(  k\left(  \theta
,t\right)  \right)  d\theta.
\end{equation}
In particular, when $\Phi\left(  k\right)  =F\left(  k\right)  -\lambda\left(
t\right)  \ $is from the LP flow (\ref{eqn1.3L})\ or from the
AP\ flow\ (\ref{eqn1.3A}), the right hand sides of (\ref{in-1}) and
(\ref{in-2}) will become zero, which gives a useful estimate. More precisely,
we have:

\begin{lemma}
(\emph{Entropy estimate.}) Under the\ LP flow\ (\ref{eqn1.3L}),\ we have%
\begin{equation}
\text{if}\ 0<\alpha<1,\ \ \int_{0}^{2\pi}k^{\alpha-1}\left(  \theta,t\right)
d\theta\ \ \text{is\ increasing in time }t\in\lbrack0,\infty), \label{40}%
\end{equation}
and%
\begin{equation}
\text{if}\ \alpha>1,\ \ \int_{0}^{2\pi}k^{\alpha-1}\left(  \theta,t\right)
d\theta\ \ \text{is\ decreasing in time }t\in\lbrack0,\infty). \label{50}%
\end{equation}
Also under the\ AP flow\ (\ref{eqn1.3A}),\ we have
\begin{equation}
\text{if}\ 0<\alpha<1,\ \ L^{\alpha-1}\left(  t\right)  \int_{0}^{2\pi
}k^{\alpha-1}\left(  \theta,t\right)  d\theta\ \ \text{is\ increasing in time
}t\in\lbrack0,\infty), \label{20}%
\end{equation}
and%
\begin{equation}
\text{if}\ \alpha=1,\ \ \int_{0}^{2\pi}\log\left(  k\left(  \theta,t\right)
L\left(  t\right)  \right)  d\theta\ \ \text{is\ decreasing in time }%
t\in\lbrack0,\infty), \label{30}%
\end{equation}
and%
\begin{equation}
\text{if}\ \alpha>1,\ \ L^{\alpha-1}\left(  t\right)  \int_{0}^{2\pi}%
k^{\alpha-1}\left(  \theta,t\right)  d\theta\ \ \text{is\ decreasing in time
}t\in\lbrack0,\infty). \label{10}%
\end{equation}

\end{lemma}%

\proof
We use (\ref{dkdt}) and apply (\ref{in-1}) in Lemma \ref{lem-Mink} to the LP
flow\ (\ref{eqn1.3L}) to get%
\begin{align*}
&  2\pi\int_{0}^{2\pi}\frac{k^{\alpha}\left(  \theta,t\right)  -\lambda\left(
t\right)  }{k^{2}\left(  \theta,t\right)  }\frac{\partial k}{\partial
t}\left(  \theta,t\right)  d\theta\\
&  \leq\left(  \int_{0}^{2\pi}\left(  k^{\alpha}\left(  \theta,t\right)
-\lambda\left(  t\right)  \right)  d\theta\right)  ^{2}=\left(  \frac{dL}%
{dt}\left(  t\right)  \right)  ^{2}=0,
\end{align*}
where$\ \lambda\left(  t\right)  =\left(  2\pi\right)  ^{-1}\int_{0}^{2\pi
}k^{\alpha}\left(  \theta,t\right)  .\ $In such a\ case, we will obtain
\begin{equation}
\frac{d}{dt}\left(  \frac{1}{\alpha-1}\int_{0}^{2\pi}k^{\alpha-1}\left(
\theta,t\right)  d\theta\right)  \leq0,
\end{equation}
which gives (\ref{40})\ and (\ref{50}).

On the other hand, apply (\ref{in-2})\ in Lemma \ref{lem-Mink} to the AP
flow\ (\ref{eqn1.3A}) to get%
\begin{align*}
&  2A\left(  t\right)  \int_{0}^{2\pi}\frac{k^{\alpha}\left(  \theta,t\right)
-\lambda\left(  t\right)  }{k^{2}\left(  \theta,t\right)  }\frac{\partial
k}{\partial t}\left(  \theta,t\right)  d\theta\\
&  \leq\left(  \int_{0}^{2\pi}\frac{k^{\alpha}\left(  \theta,t\right)
-\lambda\left(  t\right)  }{k\left(  \theta,t\right)  }d\theta\right)
^{2}=\left(  \frac{dA}{dt}\left(  t\right)  \right)  ^{2}=0,
\end{align*}
where$\ \lambda\left(  t\right)  =L^{-1}\left(  t\right)  \int_{0}^{2\pi
}k^{\alpha-1}\left(  \theta,t\right)  d\theta.\ $The above is the same as%
\[
\int_{0}^{2\pi}k^{\alpha-2}\left(  \theta,t\right)  \frac{\partial k}{\partial
t}\left(  \theta,t\right)  d\theta+\frac{{L}^{\prime}\left(  t\right)  }%
{{L}\left(  t\right)  }\int_{0}^{2\pi}k^{\alpha-1}\left(  \theta,t\right)
d\theta\leq0,\ \ \ L\left(  t\right)  =\int_{0}^{2\pi}\frac{1}{k\left(
\theta,t\right)  }d\theta.
\]
For $\alpha=1,\ $we get (\ref{30}).\ For $\alpha>0,\ \alpha\neq1,\ $we get%
\[
\frac{1}{\alpha-1}\frac{d}{dt}\log\left(  \int_{0}^{2\pi}k^{\alpha-1}\left(
\theta,t\right)  d\theta\right)  +\frac{d}{dt}\log L\left(  t\right)  \leq0.
\]
The result follows.$%
\hfill
\square$

\begin{remark}
As a comparison, in Gage-Hamilton's curve shortening flow (see \cite{GH}) the
integral%
\begin{equation}
\int_{0}^{2\pi}\log\left[  k\left(  \theta,t\right)  \sqrt{\frac{A\left(
t\right)  }{\pi}}\right]  d\theta
\end{equation}
is decreasing in time, which is called \emph{entropy estimate}.\ Now, in the
area-preserving\ curve shortening flow\ (see \cite{Ga2}) the integral%
\begin{equation}
\int_{0}^{2\pi}\log\left[  k\left(  \theta,t\right)  L\left(  t\right)
\right]  d\theta
\end{equation}
is decreasing in time. \ 
\end{remark}

\ \ 

\ \ \ \ \ \ 

\ \ \ \ \ \ 

\ \textbf{Acknowledgment.\ \ }\ \ We thank the referee for his careful reading
of our paper and useful comments and suggestions.\ We are very grateful to
Professor Ben Andrews of the Australian National University.\ Without
consulting him so often, it is unlikely for us to finish this paper. The
first\ author is supported by the National Science Council\ of Taiwan with
grant number 102-2115-M-007-012-MY3.\ The second author is supported by the
National Natural Science Foundation of\ China\ 11101078, 11171064,\ and the
Natural Science Foundation of\ Jiangsu Province BK20130596.

\ \

\ \ \ \ \ 

\ \ 

\ 

Dong-Ho Tsai

Department of Mathematics

National Tsing Hua University

Hsinchu 300,\ TAIWAN

E-mail:\ \textit{dhtsai@math.nthu.edu.tw}

\ \ \ \ 

\ \ \ 

Xiao-Liu Wang

Department of Mathematics

Southeast University

Nanjing 210096, PR CHINA

E-mail:\ \textit{xlwang@seu.edu.cn}

\end{document}